\documentclass[a4,10pt]{amsart}
\usepackage{amsmath}
\usepackage{amssymb}

\begin{document}

\numberwithin{equation}{section}

\newtheorem{theorem}{Theorem}[section]
\newtheorem{assertion}[theorem]{Assertion}
\newtheorem{claim}[theorem]{Claim}
\newtheorem{conjecture}[theorem]{Conjecture}
\newtheorem{corollary}[theorem]{Corollary}
\newtheorem{defn}[theorem]{Definition}
\newtheorem{example}[theorem]{Example}
\newtheorem{figger}[theorem]{Figure}
\newtheorem{lemma}[theorem]{Lemma}
\newtheorem{prop}[theorem]{Proposition}
\newtheorem{remark}[theorem]{Remark}
\newtheorem{assumption}[theorem]{Assumption}

\title
  { The Kahler Metrics of constant scalar curvature on the Complex Torus}

\author{ Bohui Chen$^1$}
\address{ Yangtze Center of Mathematics, Department of Mathematics, Sichuan University\\
        Chengdu}
\email{bohui\@cs.wisc.edu}

 \author  { An-Min Li$^{2,*}$}
\address{
 Yangtze Center of Mathematics, Department of Mathematics, Sichuan University\\
        Chengdu}
\email{math$\_$li\@yahoo.com.cn} \footnote[1]{partially supported
by NKBRPC(2006CB805905), NSFC 10631050 }
 \footnote[2]{partially supported by
NKBRPC(2006CB805905), NSFC 10631050  and
RFDP(20060610004)}
\renewcommand{\thefootnote}{\fnsymbol{footnote}}
\setcounter{footnote}{0} \footnote{Corresponding author.}

 \author  { Li Sheng }
\address{
 Department of Mathematics, Sichuan University\\
        Chengdu}
        \email{l$\_$sheng\@yahoo.cn}
\noindent
\def \mc {\mathcal}
 \def \J{{\cal J}}
              \def \Map{Map(S^2, V)}
              \def \M{{\cal M}}
              \def \A{{\cal A}}
              \def \B{{\cal B}}
              \def \C{{\bf C}}
              \def \Z{{\bf Z}}
              \def \R{{\bf R}}
              \def \P{{\bf P}}
              \def \I{{\bf I}}
              \def \N{{\bf N}}
              \def \T{{\bf T}}
              \def \O{{\cal O}}
              \def \Q{{\bf Q}}
              \def \D{{\bf D}}
              \def \H{{\bf H}}
              \def \s{{\mathcal S}}
              \def \e{{\bf E}}
              \def \k{{\bf k}}
              \def \U{{\cal U}}
              \def \E{{\cal E}}
              \def \F{{\cal F}}
              \def \L{{\cal L}}
              \def \K{{\cal K}}
              \def \G{{\bf G}}

\def \real{\mathbb{R}}

\maketitle

{\bf Abstract. We study the Dirichlet problem of the Abreu equation. The solutions
provide the Kahler metrics of constant scalar curvature on the complex torus.}

\vskip 0.1in \noindent \vskip 0.1in \noindent

One of the central problem in complex geometry is to find certain
canonical metrics within a given Kahler class. As examples, the
extremal metrics, introduced by E. Calabi, has been studied
intensively in the past 20 years. Most extremal metrics are Kahler
metrics of constant scalar curvature. There are three aspects of the
problem: sufficient conditions of existence, necessary condition of
existence and uniqueness.  The necessary conditions for the
existence are conjectured to be related to certain stabilities. For
example, it was first by Tian that gave an analytic "stability"
condition which is equivalent to the existence of a Kahler-Einstein
metric (\cite{T}). There are many works on this aspects (\cite{T},
\cite{D1},\cite{D2},\cite{CT}).  The uniqueness aspect is completed
by Mabuchi (\cite{M})  in the algebraic case and Chen-Tian
(\cite{CT}) in general, in the sense that the extremal metric is
unique up to the action of holomorphic automorphisms.

On the other hand, there has been not much progress on the
existence of extremal metrics or Kahler metrics of constant scalar
curvature. One reason is that the equation is highly nonlinear and
of 4th order. Our project is to understand this problem on toric
varieties following the works of Abreu and Donaldson
(\cite{Ab1},\cite{D3},\cite{D4}). When studying the equation for
Kahler metrics with prescribed scalar curvature on {\em toric
varieties}, one can reduce the equation of complex variables to a
real equation    on a polytope in $\mathbb{R}^n$. In \cite{Ab1},
using Guillemin's method(\cite{G}),
 Abreu formulates this equation \eqref{eqn_1.1} which is called
the Abreu equation now. Since this is a 4th order equation, the progress
on this equation is slow. One of the main result is given by Donaldson.
He gives the interior estimates of the equation when $n=2$ (\cite{D4}).

In this paper, we study a PDE problem: the Dirichlet problem of the
Abreu equation on strictly convex domain with degenerated boundary
conditions. We show the interior regularity of the  equation.  As a
corollary, we have constructed abundant Kahler metrics of constant
scalar curvature
 on complex torus $(\mathbb{C}^\ast)^n$,
where $\mathbb{C}^\ast=\mathbb{C}-\{0\}$. (cf. Remark \ref{rmk_1.2}
and Corollary \ref{cor_1.1}.) The graphs of solution we construct are
 Euclidean complete, however
presumably the associated  metrics are
not complete on the complex torus.

The paper is organized as following: in the introduction section
we review the equations,
formulate our problem and state the main theorem in this paper; the estimates
of determinant are given in 2nd section, in particular,
Lemma \ref{lemma_2.3} is the core lemma of the paper; the rest of paper is
devoted to the proof of the main theorem.

\section{Introduction}\label{sect_1}
\vskip 0.1in
\noindent Given a
bounded convex domain $\Omega \subset R^{n}$, we study the
Abreu equation in this paper
$$\s(u) = K$$
where $K(\xi), \xi=(\xi_1,...,\xi_n)$ is some given smooth function defined
on an open subset of $\mathbb{R}^{n}$ containing $\overline {\Omega}$
and $\s(u)$ denotes the expression
$$
\s(u) = - \sum \frac{\partial u^{ij}}{\partial \xi_i \partial
\xi_j}.
$$
The Abreu equation appears in the study of the differential geometry
of toric varieties ( see [A], [D-1],[D-2]), where $K$ is the scalar
curvature of the Kahler metric. The Kahler metric is extremal in the
Calabi sense if and only if $K$ is an affine function in $\xi$. The
Abreu equation can be written as (see [D-4] section  2.1)
\begin{equation}\label{eqn_1.1}
\sum_{i,j=1}^n U^{ij} w_{ij} = - K\;\;\;in\;\Omega
\end{equation}
where $(U^{ij})$ is the cofactor matrix of the Hessian matrix $D^2
u$ of the convex function $u$ and
$
w =\det(u_{ij})^{-1}.$

For any smooth and strictly convex function
$u$ on $\Omega$ we consider the normal map
$$L_u : \Omega\rightarrow L_u(\Omega)=:\Omega^\ast \subset \real^n$$
defined by
$$(\xi_1,...,\xi_n)\mapsto (x_1,...,x_n)$$
$$x_i = \frac{\partial
u}{\partial \xi_i}. $$ Then $L_u$ is a diffeomorphism. Define a
function $f(x)$ on $\Omega^\ast$ by
\begin{equation}\label{eqn_1.3}
f(x) = \sum_{k=1}^n \xi_k\frac{\partial u}{\partial \xi_k} - u(\xi).
\end{equation}
$f$ is called the {\em Legendre transformation} of $u$.
In terms of $x_i$ and $f(x)$ the Abreu equation can be written as
\begin{equation}\label{eqn_1.4}
K = -\sum_{i,j=1}^n f^{ij}\frac{ \partial^2}{\partial x_i\partial
x_j}\log(\det(f_{kl})),
\end{equation}
where $(f^{ij})$ is the inverse of the Hessian $(f_{ij})$.

 Denote
$$
 M^\ast = \left\{(x,f(x))|x\in \Omega^\ast\right\}.
$$
to be the graph of $f$ over $\Omega^\ast$. If $|\nabla
u||_{\partial \Omega} = \infty $ , then $f(x)$ is defined on whole
$\mathbb{R}^n$, i.e., $ M^\ast$ is Euclidean complete.

The main result of this paper is the following. The proof is given
in \S\ref{sect_3}.
\begin{theorem}\label{thm_1.1}
 Let
$\Omega\subset \mathbb{R}^{n}$
be a bounded domain with smooth and strictly
convex boundary, $K$ be a smooth function defined on an open subset
of $\mathbb{R}^{n}$ containing $\overline {\Omega}$ such that $K \geq k_o
>0$ for some constant $k_o>0$. Given a smooth and strictly convex
function function $\varphi$ defined on an open subset of $\mathbb{R}^{n}$
containing $\overline {\Omega}$, then there is a function $u$ such
that
\begin{itemize}
\item $u$ is smooth and strictly convex in $\Omega$;
\item on $\partial\Omega$
 $$
 u = \varphi,\;\;\;\;
|\nabla u| = \infty, \;\;\;w = 0;
$$
\item $u$ solves the Abreu equation
$$\s(u) = K\;\;\;in\;\Omega .$$
\end{itemize}
\end{theorem}

\begin{remark}\label{rmk_1.2}
Let $\Omega$ be a convex domain and $u$ be a smooth and strictly
convex function that solves the Abreu equation in $\Omega$.
Suppose that $f$ is the legendre transformation of $u$ and
$L_u(\Omega) =\mathbb{R}^n$. Then $f$, as a potential function of
the Kahler matric on the complex torus $\mathbb{C}^n
/2\pi\sqrt{-1}\mathbb{Z}^n, $ gives a $T^n$-invariant  Kahler
metric of scaler curvature $K$.
\end{remark}
Hence we have
\begin{corollary}\label{cor_1.1}
Let $u$ be a solution of the Abreu equation given in Theorem
\ref{thm_1.1}. Then the Legendre transform
function $f$ of $u$ yields a
 Kahler metric of  scaler curvature $K$. In particular,
if $K$ is constant, we construct a $T^n$-invariant  Kahler metric
of constant curvature $K$ on the complex torus $\mathbb{C}^n
/2\pi\sqrt{-1}\mathbb{Z}^n$.
\end{corollary}
{\bf Proof. }Since $|\nabla u|=\infty$ at boundary, $f$ defines on
the whole $\mathbb{R}^n$. Therefore the claim follows from  Remark
\ref{rmk_1.2}. $\Box$

\section{Estimates for determinant}\label{sect_2}
In this section we derive some estimates of determinant
$\det(D^2u)$. The following two lemmas can be found in [D-4],
[T-W-1].
\begin{lemma}\label{lemma_2.1}
Suppose that $u$ is a smooth and strictly convex function in $\Omega$
satisfying the Abreu equation $\s(u) = K$. If $\det(u_{ij}) > d_1$
near $\partial \Omega$, then
$$\det (u_{ij})\geq d_1$$
everywhere in $\Omega$, where
$$d_1 = \left(\frac{4Max_{\Omega}\{K\}
 diam(\Omega)^2}{n}\right)^{-n}.$$
\end{lemma}
For any $p\in \Omega$ we say that $u$  is normalized at $p$ if
$$u \geq 0,\;\;\;\;u(p) = 0,\;\;\;\frac{\partial u}{\partial \xi_k}(p) = 0 \;\;\;\;\forall \;k=1,...,n.$$
On the other hand, suppose that $u$ is not normalized at $p$. Let
$$
\xi_{n+1}=a\cdot (\xi-p)+ b
$$
be the support hyperplane of $u$ at $p$. Set
$$
\tilde u= u-(a\cdot(\xi-p)+b).
$$
Then $\tilde u$ is normalized at $p$. We call $\tilde u$ is a {\em normalization}
of $u$ with respect to $p$.

Suppose that $u$ is normalized at $p$. For any positive number $b$
we denote
$$S_u(p,b) = \left\{ \xi\in\Omega | u(\xi)< b \right \},$$
$$\bar{S}_u(p,b) = \left\{ \xi\in\Omega | u(\xi)\leq b \right \}.$$
\begin{lemma}\label{lemma_2.2}
Suppose that $u$ is a smooth and strictly convex  function defined
in $\Omega$ with $\s(u) = K$.
Suppose that $u$ is normalized at $p$ in $\Omega$. If the section
$$\bar{S}_u(p,C)=\{\xi\in \Omega| u(\xi)\leq C\}$$
is compact and if there is a constant $b>0$ such that
$$\sum _{k=1}^n x_k^2 \leq b$$
on $\bar{S}_u(p,C)$, then there is a constant $d_2>0$ depending on
$n$, $C$ and $b$ such that the estimate
$$\det (u_{ij})\leq d_2$$ holds in $S_u(p,C/2)$.
\end{lemma}

Let $\Omega \subset \mathbb{R}^n$ be a bounded, normalized   convex
domain. Thus
\begin{equation}\label{eqn_0.1}
n^{-\frac{3}{2}}D_1(0) \subset  {\Omega} \subset
D_1(0).\end{equation}

 Denote by $\mathcal{F}(\Omega,C)$ the class of convex
functions defined on $\Omega$ such that
$$ \inf_{\Omega} {u } = 0,\;\;\;
u= C\;\;on\;\;\partial \Omega.$$
\begin{lemma}   Let $\Omega_k$ be a sequence of smooth and normalized convex
domains, $u_k\in \mc F (\Omega_k,C)$. Then there are
constants $d>1$, $b>0$ independent of $k$ such that
$$\frac{\sum_i(\frac {\partial u_k}{\partial \xi_i})^2 }
{(d+f_k)^2}\leq b, \ \ \ \ k=1,2,\dots \quad\hbox{on}\ \
\bar{\Omega}_k.$$ \end{lemma}

{\bf Proof.} We may suppose by taking subsequence that $\Omega_k$
converges to a convex domain $\Omega $ and $u_k$ converges to a
convex function $u_\infty $, locally uniformly in $\Omega$.
Obviously, we have the uniform estimate
\begin{equation}\sum \left(\frac{\partial u_k}{\partial \xi_i}\right)^2(0)\leq
4n^3.\end{equation} For any $k$, let
\begin{equation}\tilde{u}_k = u_k - \sum \frac{\partial
u_k}{\partial \xi_i}(0)\xi_i - u_k(0).\end{equation}
 Then
$$\tilde{u}_k(0) = 0,\;\;\;\;\tilde{u}_k(\xi)\geq 0,\;\;\;
\tilde{u}_k|_{\partial \Omega_k}\leq C_0,$$ where $C_0$ is a
constant depending only on $n$. As $B(0,n^{-\frac{3}{2}})
\subset\Omega_k$, we have
\begin{equation}\frac{\mid\nabla \tilde{u}_k\mid^2}{(1+\tilde{f}_k)^2}\leq
\mid\nabla \tilde{u}_k\mid^2\leq
\frac{C_0^2}{dist(B(0,2^{-1}n^{-\frac{3}{2}}),\partial\Omega_k)^2}\leq
4n^{3}C^2_0 \ \ \ \ \ \end{equation} on
$B(0,2^{-1}n^{-\frac{3}{2}})$, where
 $\tilde{f}_k$ is the
Legendre transformation of $\tilde{u}_k$ relative to $0$.
We discuss three cases.

{\bf Case 1.}   $p\in B(0,2^{-1}n^{-\frac{3}{2}})$. Following from
(2.2), (2.3) and (2.4) we have
$$\frac{\mid\nabla u_k\mid^2}{(1+f_k)^2}\leq
4n^{3}(C^2_0+1).
$$
{\bf Case 2.} $p\in\partial{\Omega}_k$, we
may suppose that $p =(\xi_1,0,\dots,0)$ with $\xi_1> 0$ by an
orthonormal transformation. Then, at $p$,
$$C_0+\tilde{f}_k\geq\tilde{u}_k+\tilde{f}_k=
\frac{\partial \tilde{u}_k}{\partial\xi_1}\xi_1.$$ It follows
that
$$\frac{\left(\frac {\partial \tilde{u}_k}{\partial \xi_1}\right)^2 }
{(C_0+\tilde{f}_k)^2}<\frac{1}{\xi_1^2}<4n^3.$$ Therefore
there exist constants $\tilde{d}>1,$ $\tilde{b}>0$ depending only
on $n$ such that
\begin{equation}\frac{\left(\frac {\partial \tilde{u}_k}{\partial r}\right)^2 }
{(\tilde{d}+\tilde{f}_k)^2}< \tilde{b},\end{equation} where
$\frac{\partial}{\partial r}$ denotes the radial derivative.
 Note that
\begin{equation} \frac{\partial\tilde{u}_k}{\partial \xi_i}= \frac{\partial
u_k}{\partial \xi_i}- \frac{\partial{u}_k}{\partial
\xi_i}(0),\ \ \ \tilde{f}_k= f_k + u_k(0).
 \end{equation}
It follows from (2.3) and (2.4) that
$$\left(\frac {\partial u_k}{\partial r}\right)^2 \leq
2\left(\frac {\partial \tilde{u}_k}{\partial r}\right)^2 +
8n^3.$$ Then
$$\frac{\left(\frac {\partial u_k}{\partial r}\right)^2 }
{(d'+ f_k)^2}< b',\leqno(5.8)$$ for some constants $d'>1$,
$b'>0$ independent of $k$. Note that
$$|\nabla u_k(p)|=\frac{1}{\cos\alpha_k}\left|\frac
 {\partial u_k}{\partial r}(p)\right|,\leqno(5.9)$$
where $\alpha_k$ is the angle between vectors $\nabla u_k(p)$
and $\frac {\partial u_k}{\partial r }(p).$ Since $u_k= 1
$ on $\partial \Omega_k$, $\nabla u_k(p)$ is perpendicular to
the boundary of $\Omega_k$ at any $p\in \partial \Omega_k$. As
$\Omega$ is convex and $0\in \Omega$, it follows that
$\frac{1}{\cos\alpha_k}$ have a uniform upper bound.

{\bf Case 3.} $p\in {\Omega}_k^o\backslash
B(0,2^{-1}n^{-\frac{3}{2}}).$ Let $F_k=\frac{\sum
x_k^2}{(d+f_k)^2}.$ We can assume that
  \begin{equation}
\label{eqn_2.7.b}
  \;\;\;\; \max_{\partial\Omega_k\cup \bar B(0,2^{-1}n^{-\frac{3}{2}})}F_k <\max_{\bar {\Omega}_k} F_k. \end{equation}
In fact, if \eqref{eqn_2.7.b} is not true, then the lemma follows
from a direct calculation. Let $p^\ast_k\in {\Omega}_k^o\backslash
B(0,2^{-1}n^{-\frac{3}{2}})$ be the point such that
$$F_k(p^\ast_k)=\max_{\Omega_k}F_k.$$ Then, at $p^\ast_k,$  $\frac{\partial}{\partial x_i}F_k=0.$ Thus,
\begin{equation}
\frac{x_i}{\sum x_j^2}=\frac{\xi_i}{d+f_k}.
\end{equation}
where $d=d'+C.$ Obviously, $d+f_k\geq 1.$ Then
$$\frac{\max_{i=1}^n|x_i|}{\sum x_j^2}=\frac{\max|\xi_i|}{d+f_k}\geq \frac{2^{-n}n^{-2}}{d+f_k}.$$ On the other hand,
$$\frac{\max_{i=1}^n|x_i|}{\sum x_j^2}\leq \frac{1}{ (\sum x_j^2 )^\frac{1}{2}}.$$
Noting that $F_k(p)\leq F_k(p^\ast_k),$ we have
$$F_k(p)\leq 4^nn^4.\;\;\;\;\; q.e.d.$$

In the following we prove a  stronger estimate than that in Lemma \ref{lemma_2.2}
which plays an important role in this paper.
\begin{lemma}\label{lemma_2.3}
Let $u$ be a smooth and strictly convex function defined in
$\Omega$ which satisfies the Abreu equation $\s(u) = K$. Suppose
that $u$ is normalized at $p$ and the section $\bar{S}_u(p,C)$ is
compact. And suppose that  there is a constant $b>0$ such that
\begin{equation}\label{eqn_2.0}
\frac{\sum x_k^2}{(1+ f)^2}\leq b
\end{equation}
on $\bar{S}_u(p,C)$. Then there is a constant $d_3
>0$ depending only on $n$,$b$ and $\frac{1}{C}$, such that
\begin{equation}\label{eqn_2.1}
\exp \left\{-\frac{4C}{C-u }\right\} \frac{\det
(u_{ij})}{(1+f)^{2n}}\leq  d_3
\end{equation}
 on $S_u(p,C)$.
\end{lemma}
{\bf Proof.}
Consider the following function
$$
F = \exp \left\{-\frac{m}{C-u } + L\right\}\frac{1}{w(d+f)^{2n}},$$
where
$$
L = \epsilon \frac{\sum x_k^2}{(1+f)^{2}}.
$$
$m$ and $\epsilon$ are positive constants to be determined later.
Clearly,
F attains its supremum at some interior point $p^*$ of $S_u(p,C)$.
We have, at $p^*$,
\begin{equation}\label{eqn_2.2}
F_{i} = 0,
\end{equation}
\begin{equation}\label{eqn_2.3}
\sum u^{ij}F_{ij} \leq 0,
\end{equation}
where we denote $F_i = \frac{\partial F}{\partial \xi_i}$, $F_{ij} =
\frac{\partial^2 F}{\partial \xi_i\partial \xi_j},f_i =
\frac{\partial f}{\partial \xi_i}$ and so on. We calculate both
expressions \eqref{eqn_2.2} and \eqref{eqn_2.3} explicitly:
\begin{equation}\label{eqn_2.4}
-\frac{m}{(C-u)^2}u_i + L_i - 2n\frac{f_i}{1+f} - \frac{w_i}{w} = 0,
\end{equation}
and
\begin{eqnarray}\label{eqn_2.5}
&&-\frac{2m}{(C-u)^3}\sum u^{ij}u_iu_j - \frac{mn}{(C-u)^2} + \sum u^{ij}{L_{ij}}
\\
&&- 2n \frac{\sum u^{ij}f_{ij}}{1+f}
 + 2n\frac{\sum u^{ij}f_if_j}{(1+f)^2} + \frac{\sum u^{ij}w_iw_j}{w^2} + K \leq 0.
 \nonumber
\end{eqnarray}
Since
\begin{eqnarray*}
f_i &=& \sum \xi_ku_{ki},\;\;\;
f_{ij} =u_{ij} + \sum \xi_ku_{kij}.
\end{eqnarray*}
Then
\begin{equation*}
- 2n \frac{\sum u^{ij}f_{ij}}{1+f} =
- \frac{2n^2}{1+f} + \frac{2n}{1+f}\sum \xi_k\frac{w_k}{w}.
\end{equation*}
By \eqref{eqn_2.4}
\begin{equation}
- 2n \frac{\sum u^{ij}f_{ij}}{1+f}
 =  - \frac{2n^2}{1+f} - \frac{2mn\sum u_i\xi_i}{(C-u)^2(1+f)} + \frac{2n}{1+f} \sum \xi_iL_i
- 4n^2 \frac{\sum \xi_if_i}{(1+f)^2}.
\end{equation}
Note that
$$\sum \xi_kf_k = \sum u_{kl}\xi_k\xi_l = \sum u^{ij}f_if_j,$$
Hence
$$\frac{2mn}{(C-u)^2(1+f)}\sum u_i\xi_i = \frac{2mn(u+f)}{(C-u)^2(1+f)}\leq \frac{2mn}{(C-u)^2}Max\{1,C\},$$
\begin{eqnarray*}\frac{2n}{1+f}
\sum \xi_iL_i &=& \frac{2n}{1+f} \sum \xi_ku_{kj}u^{ji}L_i\leq \frac{1}{2}\frac{\sum u_{kl}\xi_k\xi_l}{(1+f)^2} + 2n^2\sum u^{ij}L_iL_j
\\
& =& \frac{1}{2}\frac{\sum u^{ij}f_if_j}{(1+f)^2} + 2n^2\sum u^{ij}L_iL_j.
\end{eqnarray*}
We have
\begin{eqnarray}\label{eqn_2.7}
- 2n \frac{\sum u^{ij}f_{ij}}{1+f}
&\geq& -2n^2 - \frac{2mn}{(C-u)^2}Max\{1,C\} - 2n^2\sum u^{ij}L_iL_j\\
&& - \left(4n^2 + \frac{1}{2}\right)\frac{\sum u^{ij}f_if_j}{(1+f)^2}.
\nonumber\end{eqnarray}
Let us calculate the terms $\sum u^{ij}L_iL_j$ and $\sum u^{ij}L_{ij}$.
Since
$$L_i = \epsilon \frac{2\sum x_ku_{ki}}{(1+f)^2} - 2\epsilon \frac{f_i\sum x_k^2}{(1+f)^3},$$
then
\begin{equation}\label{eqn_2.8}
\sum u^{ij}L_iL_j \leq 8\epsilon b \frac{\epsilon \sum u_{kk}}{(1+f)^2}+ 8(\epsilon b)^2\frac{\sum u^{ij}f_if_j}{(1+f)^2},
\end{equation}
\begin{eqnarray}\label{eqn_2.9}
\sum u^{ij}L_{ij}& =& \epsilon \frac{2\sum u_{kk}
 + 2\sum x_ku^{ij}u_{ijk}}{(1+f)^2} -
 8\epsilon \frac{\sum x_kf_iu_{kj}u^{ij}}{(1+f)^3}
\\
 &&
 - 2\epsilon \frac{\sum x_k^2 \sum u^{ij}f_{ij}}{(1+f)^3} + 6\epsilon \frac{\sum x_k^2 \sum u^{ij}f_if_j}{(1+f)^4}
 \nonumber\end{eqnarray}
Applying the Schwarz inequality for each term in \eqref{eqn_2.9}:
\begin{equation*}\label{eqn_2.10}
 8\epsilon \frac{\sum x_kf_iu_{kj}u^{ij}}{(1+f)^3} \leq 16\epsilon^2\frac{\sum u_{kl}x_kx_l}{(1+f)^4} + \frac{\sum u^{ij}f_if_j}{(1+f)^2}
\leq 16\epsilon b \frac{\epsilon \sum u_{kk}}{(1+f)^2}  + \frac{\sum u^{ij}f_if_j}{(1+f)^2},
\end{equation*}
\begin{eqnarray*}\label{eqn_2.11}
2\epsilon \frac{\sum x_k^2 \sum u^{ij}f_{ij}}{(1+f)^3}&=&
 2n\epsilon\frac{\sum x_k^2}{(1+f)^3}- 2\epsilon \frac{\sum x_k^2}{(1+f)^3}
\sum \xi_k\frac{w_k}{w} \\
&\leq&
 2n\epsilon\frac{\sum x_k^2}{(1+f)^3} + \frac{1}{8n}\frac{\sum u^{ij}w_iw_j}{w^2} + 8n\left[\epsilon
\frac{\sum x_k^2}{(1+f)^2}\right]^2 \frac{\sum u_{ij}\xi_i\xi_j}{(1+f)^2}\nonumber\\
& \leq&
 2n\epsilon b + \frac{1}{8n}\frac{\sum u^{ij}w_iw_j}{w^2} + 8n (\epsilon b)^2 \frac{\sum u^{ij}f_if_j}{(1+f)^2},
\nonumber\end{eqnarray*}
\begin{eqnarray*}\label{eqn_2.12}
2\epsilon \frac{\sum x_ku^{ij}u_{ijk}}{(1+f)^2}& =&
 -2\frac{\epsilon}{(1+f)^2}\sum x_k\frac{w_k}{w}= -2\frac{\epsilon}{(1+f)^2}\sum x_lu_{li}u^{ik}\frac{w_k}{w}
\\
&\leq&
\frac{1}{8n}\frac{\sum u^{ij}w_iw_j}{w^2} + 8n\frac{\epsilon^2\sum u_{kl}x_kx_l}{(1+f)^4}\nonumber\\
&\leq& \frac{1}{8n}\frac{\sum u^{ij}w_iw_j}{w^2} + 8n\epsilon b \frac{\epsilon \sum u_{kk}}{(1+f)^2}.
\nonumber
\end{eqnarray*}
Then
\begin{eqnarray}\label{eqn_2.13}
\sum u^{ij}L_{ij} &\geq& (2-16\epsilon b - 8n\epsilon b)\frac
{\epsilon\sum u_{kk}}{(1+f)^2} - \frac{1}{4n}\frac{\sum u^{ij}w_iw_j}{w^2}
\\
&& -
\left(1+8n(\epsilon b)^2\right)\frac{\sum u^{ij}f_if_j}{(1+f)^2} - 2nb\epsilon.
\nonumber
\end{eqnarray}
Note that
\begin{equation}\label{eqn_2.14}
\frac{|\sum u^{ij}u_if_j|}{1+f} = \frac{|\sum x_k\xi_k|}{1+f} =
\frac{|u + f|}{1+f}\leq 1.
\end{equation}
By \eqref{eqn_2.4},\eqref{eqn_2.8} and \eqref{eqn_2.14} we have
\begin{eqnarray}\label{eqn_2.15}
&&\left(1-\frac{1}{4n}\right)\frac{\sum u^{ij}w_iw_j}{w^2}
\geq \left(1-\frac{1}{4n}\right)(1-\delta)\frac{m^2}{(C-u)^4}
\sum u^{ij}u_iu_j
\\
&&\;\;\;\;\;\;\;\;\;\; + \left(1-\frac{1}{4n}\right)4n^2(1-\delta)\frac{\sum u^{ij}f_if_j}{(1+f)^2} - \left(\frac{1}{\delta}-1\right)\sum u^{ij}L_iL_j
\nonumber \\
&&\;\;\;\;\;\;\;\;\;\;\geq
\left(1-\frac{1}{4n}\right)(1-\delta)\frac{m^2}{(C-u)^4}\sum u^{ij}u_iu_j
\nonumber\\
&&\;\;\;\;\;\;\;\;\;\; + \left(1-\frac{1}{4n}\right)4n^2(1-\delta)\frac{\sum u^{ij}f_if_j}{(1+f)^2}
 - \frac{8\epsilon b}{\delta}\frac{\epsilon\sum u_{kk}}{(1+f)^2}
- \frac{8(\epsilon b)^2}{\delta}\frac{\sum u^{ij}f_if_j}{(1+f)^2}
\nonumber
\end{eqnarray}
for any small positive number $\delta$. We choose  $\epsilon =
\frac{1}{8000n^2b}$, $\delta = \frac{1}{200n^2}$, $m= 4C$.
Obviously,
\begin{equation}\label{eqn_2.16}
\left(1-\frac{1}{4n}\right)\left(1-\frac{1}{24n^2}\right)\frac{m^2}{(C-u)^4} > \frac{2m}{(C-u)^3}.
\end{equation}
Inserting \eqref{eqn_2.7}, \eqref{eqn_2.12}, \eqref{eqn_2.13},
\eqref{eqn_2.15} and \eqref{eqn_2.16} into
\eqref{eqn_2.5}, we get
$$
\epsilon \frac{\sum u_{kk}}{(1+f)^2}
- \frac{3mn}{(C-u)^2}Max\{1,C\}- 3n^2 + K \leq 0.
$$
As
$$\sum u_{ii} \geq n [\det(u_{ij})]^{1/n} = n w^{-1/n}$$ we get
$$ \exp \left\{-\frac{m}{C-u } + \epsilon \frac{\sum x_k^2}{(1+
f)^{2}}\right\}\frac{1}{(1+f)^{2n}w}\leq d$$ for some constant $d
>0$ depending on $n$, $b$, $\frac{1}{C}$. Since $F$ attains its
maximum at $p^*$, (2.22) holds everywhere. $\Box$ \vskip 0.1in
\noindent We remark that $1+f$ in \eqref{eqn_2.0} can be replaced by
any $d+f$ with $d>0$. \\

Using the technique of Lemma \ref{lemma_2.3}  and that of Lemma 4.2
in \cite{L-J}, we can prove the following lemma  in the case of
$n=2.$

\begin{lemma}\label{lemma_2.3.a}
Let $u\in \mc F (\Omega,1) $  with $\mc S(u) = K$.  Suppose that
$D_r(0)\subset \Omega.$ Then there is a constant $d_3
>0$ depending only on $K_o$, $b$ and $r$, such that
\begin{equation}\label{eqn_2.1.a}
(r^2-\sum \xi_i^2)^2 \frac{\det (u_{ij})}{(d+f)^{4}}\leq  d_3.
\end{equation}
 \end{lemma}
{\bf Proof.} Consider the following function
$$
F = (r^2-\theta)^k\exp \left\{ L\right\}\frac{1}{w(d+f)^{4}},$$ in
$D_r(0),$ where
$$
\theta=\sum \xi_i^2,\;\;\;\;\; L = \epsilon \frac{\sum
x_k^2}{(d+f)^{2}}.
$$
$k$ and $\epsilon$ are positive constants to be determined later.
Clearly, F attains its supremum at some interior point $p^*$ of $
D_r(0)$. At $p^*$, we have,
\begin{equation}\label{eqn_2.4.a}
 L_i - \frac{4f_i}{d+f} -
\frac{w_i}{w}-\frac{k\theta_i}{ r^2-\theta } = 0,
\end{equation}
and
\begin{eqnarray}\label{eqn_2.5.a}
&&   \sum u^{ij}{L_{ij}} - 4 \frac{\sum u^{ij}f_{ij}}{d+f}
 + 4\frac{\sum u^{ij}f_if_j}{(d+f)^2}+ K + \frac{\sum u^{ij}w_iw_j}{w^2}
\\&&-\frac{k\sum u^{ij}\theta_{ij}}{ r^2-\theta }-\frac{ k\sum
u^{ij}\theta_{i}\theta_{j}}{\left(r^2-\theta\right)^2}\leq 0.
\nonumber
\end{eqnarray}
Then as in Lemma \ref{}  we have
\begin{eqnarray}\label{eqn_2.7.a}
\;\;\;\;\;\;\;- 4\frac{\sum u^{ij}f_{ij}}{d+f} \geq -8 - 8\sum
u^{ij}L_iL_j - \left(16 + \frac{1}{2}\right)\frac{\sum
u^{ij}f_if_j}{(d+f)^2}-\frac{8k
 \theta}{ r^2-\theta  }. \end{eqnarray}
The calculations of  the terms $\sum u^{ij}L_iL_j$ and $\sum
u^{ij}L_{ij}$ is the same as Lemma \ref{lemma_2.3}.

 By \eqref{eqn_2.4},\eqref{eqn_2.8} and \eqref{eqn_2.14} we have
\begin{eqnarray*}\label{eqn_2.15}
 \left(1-\frac{1}{8}\right)\frac{\sum u^{ij}w_iw_j}{w^2}&\geq&\left(1-\frac{1}{8}\right)16(1-\delta)\frac{\sum
u^{ij}f_if_j}{(d+f)^2}
\\
&&
 -  \frac{1}{\delta} \left(\sum
u^{ij}L_iL_j+\frac{4k^2\sum u^{ij}\xi_i\xi_j}{\left(r^2-\theta\right)^2}\right)\\
& \geq &- \frac{4k^2r^2\sum u^{ii}
}{\delta\left(r^2-\theta\right)^2} +
\left(1-\frac{1}{8}\right)16(1-\delta)\frac{\sum
u^{ij}f_if_j}{(d+f)^2}
\nonumber\\
&&
 - \frac{8\epsilon b}{\delta}\frac{\epsilon\sum u_{kk}}{(d+f)^2}
- \frac{8(\epsilon b)^2}{\delta}\frac{\sum u^{ij}f_if_j}{(d+f)^2}
\nonumber
\end{eqnarray*}
for any small positive number $\delta$. We choose  $\epsilon =
\frac{1}{2000 b}$, $\delta = \frac{1}{100}$  and $k=2$. As \eqref{}
we get
\begin{eqnarray}\label{eqn_2.17.a}
\epsilon \frac{\sum u_{kk}}{(d+f)^2} - 12 + K-\frac{2\sum
u^{ii}+8kr^2}{ r^2-\theta }-\frac{4kr^2\sum
u^{ii}}{\left(r^2-\theta\right)^2} - \frac{4k^2r^2\sum
u^{ii}}{\delta\left(r^2-\theta\right)^2} &\leq& 0. \nonumber
\end{eqnarray}
Denote
$$\Lambda=  \frac{2 }{ r^2-\theta  }+
\frac{4kr^2 }{\left(r^2-\theta\right)^2} + \frac{4k^2r^2
}{\delta\left(r^2-\theta\right)^2},\;\;\Xi=  12 + |K|+\frac{8kr^2}{
r^2-\theta } .$$ Denote by $\lambda_1,\lambda_2$ the eigenvalues of
$(u_{ij}).$
 From \eqref{eqn_2.17.a} we have
\begin{eqnarray}\;\;\;\;\;\;\;
\epsilon \frac{\lambda_1+\lambda_2}{(d+f)^2} -\Xi
 \leq  \Lambda
\left(\frac{1}{\lambda_1}+\frac{1}{\lambda_2}\right) .
\end{eqnarray} Then
\begin{eqnarray}\label{eqn_2.18.a}\;\;\;\;\;\;\;
\epsilon \frac{\lambda_1 \lambda_2}{(d+f)^2} -\Xi\sqrt{\lambda_1
\lambda_2}  \leq \Lambda
  .
\end{eqnarray}
where we used the inequality $\lambda_1+\lambda_2\geq
2\sqrt{\lambda_1\lambda_2}$.

Multiplying $e^L( r^2-\sum \xi_i^2)^2(d+f)^{-2}$ on both sides of
\eqref{eqn_2.18.a} and applying Schwarz's inequality
  we get
$$ F\leq
d$$ for some constant $d
>0$ depending on   $b$, $r$ and $K_o$. Since $F$ attains its
maximum at $p^*$, \eqref{eqn_2.1.a} holds everywhere in $ D_r(0)$. $\Box$\\

By the same calculation of Lemma \ref{lemma_2.3.a} we can obtain
\begin{corollary}\label{lemma_2.7}
Let $\Delta\subset \mathbb{R}^{2}$ be a Delzant polytope.  Suppose
that $u\in C^{\infty}(\Delta^o)$ satisfies $|S(u)|\leq K_o$. And
suppose that there is a constant $b>0$ such that
\begin{equation}\label{eqn_2.0.a}
\frac{\sum x_k^2}{(d+ f)^2}\leq b.
\end{equation}  Then
the following estimate holds
$$\frac{\det(u_{ij})}{(d+f)^4} \leq \frac{b_0}{ dist(\xi, \partial
\Delta)^{4}}$$ for some constant $b_0>0$ depending only on $K_o$.
\end{corollary}

The following lemma is useful for the condition \eqref{eqn_2.0}.
\begin{lemma}\label{lemma_2.35}
Let $u$ be a  smooth and strictly convex function on $\Omega$
and $f$ be its Legendre transform. Suppose that $u$
is normalized at $p$ and satisfies
$$
\frac{\sum x_k^2}{(d+ f)^2}\leq b,
$$
for some constants  $d,b>0$.
For any $\tilde p$, let $\tilde u$ be the normalization of $u$ with respect to
$\tilde p$. Let $\tilde f$ be the Legendre transformation of $\tilde u$.
Then there exist constants $d'$ and $b'$ such that
$$
\frac{\sum \tilde x_k^2}{(d'+ \tilde f)^2}\leq b'.
$$
Here $\tilde x=\partial \tilde u/\partial \xi$.
\end{lemma}
{\bf Proof. }Suppose  that the support plane at $(\tilde p, u(\tilde p))$
is
$$
\xi_{n+1}=a\cdot(\xi-\tilde p) +b.
$$
Then $\tilde u= u-a\cdot(\xi-\tilde p)-b$. By direct computations, we know
that
$$
\tilde x(\xi)=x -a, \;\;\; \tilde f(\tilde x(\xi))=f(x(\xi))
-a\cdot\tilde p +b.
$$
Hence
\begin{eqnarray*}
\frac{\sum  \tilde x_k^2}{(d'+\tilde f)^2}
&=&\frac{|\tilde x|^2}{(d+ f)^2 }
\cdot \frac{(d+ f)^2}{(d'+ \tilde f)^2}\\
&\leq& \frac{2(|x|^2+|a|^2)}{(d+ f)^2}
\cdot\frac{(d+  f)^2}{(d'+ f-a\cdot\tilde p +b)^2}\\
&\leq & C(|a|)b=b'.
\end{eqnarray*}
Here, $C(|a|)$ is a constant depending on $|a|$. $d'$ is chosen so that
$d'-a\cdot \tilde p+b> d $. $\Box$.

\vskip 0.1in \noindent Let $\Delta\in R^n$ be a Delzant polytope (
for the definition of Delzant polytope please see [?]). Suppose that
$\Delta$ has $v$ vertices and $d$ faces of $(n-1)$-dimension. Denote
by $v_A$ the number of vertices in the $n-1$-dimensional face
$\ell_A = 0$, and $L= \frac{min \{v_1,...,v_d\}}{v}$. Suppose that
$$ \ell_A = \sum a_{Aj}\xi_j-\lambda_A,$$
where $1 \leq A,B,C,...\leq d$,  $1\leq i,j,k,...\leq n$.
Let $(y_1, y_2,...,y_d)$ be an affine  coordinate system in $\mathbb{R}^{d}$.
 \begin{lemma}\label{lemma_2.4}
 Set
$$P = \{(y_1,...,y_d)| y_i >0\;\; \forall \; i\}.$$
 Let $\alpha < \frac{1}{2d}$ be a positive constant, $g(y_1,...,y_d)$ be a function
 defined in $P$ given by
$$g = - \left(y_1y_2...y_d \right)^{\alpha}.$$
Then $g$ is a smooth and strictly convex function.
 \end{lemma}
{\bf Proof. } By a direct calculation we have
$$g_{AA} = - \alpha(1-\alpha)\frac{g}{y^2_{A}},\;\; A\; =\; 1,\;...\; d,$$
$$g_{AB} = \alpha^2 \frac{g}{y_Ay_B}, \;\; for\; A\ne B.$$
We claim that the matrix $\left[\left(\frac{1}{2}\delta_{AB} - \alpha \right)\frac{1}{y_Ay_B}\right]$
is positive definite. To prove this we should calculate its
all principle minors. Let $\{j_1,...,j_k\}\subset \{1,...,d\}$
such that $j_1<j_2<...<j_k$. A direct calculation gives us
\begin{equation}\label{eqn_2.17}
det \left[\left(\frac{1}{2}\delta_{j_kj_l} - \alpha \right)\frac{1}{y_{j_k}y_{j_l}}\right] =
\frac{1}{(y_{j_1}y_{j_2}...y_{j_k})^2}\frac{1}{2^k}\left(1-2k\alpha
\right) > 0.
\end{equation}
 The claim is proved. It follows that
\begin{equation}\label{eqn_2.18}
\sum _{1\leq A,B\leq d}g_{AB}h_A h_B \geq \frac{-\alpha g}{2}\sum_{1\leq A\leq d}
 \frac{h_A^2}{y_A^2}
\end{equation}
for any $(h_1,...,h_d)\in \mathbb{R}^d$. The lemma  is proved. $\Box$.
\vskip 0.1in
\noindent

Let $i: \mathbb{R}^n \rightarrow\mathbb{ R}^d$ be given by
$$y_A = \sum _{i=1}^na_{Ai}\xi_i-\lambda_A.$$
Put $\hat{u}(\xi_1,...,\xi_n) = i^*g(y_1,...,y_d), \;\;\xi \in \Delta.$

\begin{lemma}\label{lemma_2.5}
Suppose that $\alpha < min \{ \frac{1}{2d}, \frac{2L}{n}\}$. Then $\hat{u}$ is a smooth and strictly
convex function defined on $\Delta^o$, and there is a constant $c>0$ such that
$$
\det(\hat{u}_{ij}) \geq c \frac{1}{(\ell_1...\ell_d)^{2L-n\alpha}} .$$
\end{lemma}{\bf Proof. } By a linear transformation in
$GL(n,\mathbb{Z}^n)$,
we may assume that $\ell_i = \xi_i,1\leq i\leq n$.
 Then by \eqref{eqn_2.18}
\begin{equation}\label{eqn_2.19.a}
\sum \hat{u}_{ij}b_ib_j \geq \sum \frac{-\alpha
g}{2}\left[\sum_{i=1}^n \frac{b_i^2}{\xi_i^2} + \sum_{A=n+1}^d
\frac{(\sum a_{Aj}b_j)^2}{\ell_A^2}\right] \geq \sum \frac{-\alpha
g}{2}\sum_{i=1}^n \frac{b_i^2}{\xi_i^2} \end{equation} for any
$(b_1,...,b_n)\in R^n$. It follows that $\hat{u}$ is strictly convex
in $\Delta$. By \eqref{eqn_2.19.a} we have
\begin{equation}\label{eqn_2.19}
det(\hat{u}_{ij}) \geq \frac{\alpha^n}{2^n}\frac{|g|^n}{(\xi_1...\xi_n)^2}.
\end{equation}
For each vertex we have the inequality \eqref{eqn_2.19}. So
$$det(\hat{u}_{ij})^v \geq \frac{\alpha^{nv}}{2^{nv}}
 \frac{(\ell_1...\ell_d)^{\alpha nv}}{(\ell_1^{v_1}...\ell_d^{v_d})^2}.$$
Therefore there is a constant $c>0$ such that
$$det(\hat{u}_{ij}) \geq c \frac{1}{(\ell_1...\ell_d)^{2L-n\alpha}}.$$
The Lemma \ref{lemma_2.5}  is proved. $\Box$

\begin{lemma}\label{lemma_2.6}
Let $\alpha < min \{ \frac{1}{2d}, \frac{2L}{n}\}$. Suppose that $u\in C^{\infty}(\Delta)$ satisfies
the Abreu equation \eqref{eqn_1.1}. Then the following estimate holds
$$det(u_{ij}) \geq \frac{1}{b_1(\ell_1...\ell_d )^{\alpha}}$$
for some constant $b_1>0$.
\end{lemma}
 {\bf Proof }
From the Abreu's equation \eqref{eqn_1.1} and Lemma \ref{lemma_2.5}
 we can find a
constant $b_1>0$ such that
$$\sum_{1\leq i,j\leq n} U^{ij}\left( w - b_1(\ell_1...\ell_d)^{\alpha}\right)_{ij} >0.$$
Since $w = 0$ and $\hat{u} = 0$ on $\partial \Delta$, by the maximum
principle the lemma  follows. $\Box$
\begin{lemma}\label{lemma_2.7}
Let $\Omega\subset \mathbb{R}^{n}$ be a bounded domain with
smooth and strictly convex boundary. Suppose that $u\in
C^{\infty}(\Delta^o)$ satisfies the Abreu equation. Then the
following estimate holds
$$\det(u_{ij}) \geq \frac{1}{b_1 dist(\xi, \partial \Omega)^{\alpha}}$$
for some constant $b_1>0$ and $\alpha < \frac{1}{2(n+1)}$.
\end{lemma}
{\bf Proof.} For any point $\bar{\xi}\in \partial
\Omega $ we choose coordinates such that $\Omega \subset D$, where
$$D:= \{\xi | \xi_1\geq 0,..., \xi_n \geq 0, \xi_1 + \cdots + \xi_ n
\leq a \}$$ and $\bar{\xi} = (0, a_1,...,a_n)$.
 As before we may choose $d>0$ such that
$$\sum_{1\leq i,j\leq n} U^{ij}\left( w - d(\xi_1...\xi_n \ell)^{\alpha}\right)_{ij} >0,$$
where $\ell = a -\sum_{1\leq k\leq n} \xi_k$, $\alpha = \frac{1}{2(n+1)}$. Since
$$w -  d(\xi_1...\xi_n \ell)^{\alpha} \leq 0 \;\;\;on \;\partial
\Omega,$$
$$w -  d(\xi_1...\xi_n \ell)^{\alpha} = 0\;\;\; at\; \bar{\xi},$$
we have
$$w \leq  d(\xi_1...\xi_n \ell)^{\alpha}.$$
It follows that
$$det(u_{ij}) \geq \frac{1}{d_1 \xi_1^{\alpha}}$$
for some constant $d_1 >0$. Since $\bar{\xi}$ is arbitrary, by
compactness the lemma  is proved. $\Box$

\begin{remark}\label{rmk_2.1}
In the proof of Theorem \ref{thm_1.1} we will consider the
 perturbational Abreu equation
\begin{equation}\label{eqn_2.20}
\sum U^{ij} w_{ij} = - K \;\;\;
det(u_{kl})= w^{-1+ \theta}
\end{equation}
where $\theta $
is a very small positive number. It is easy to see that the above
Lemma \ref{lemma_2.1} - Lemma \ref{lemma_2.7} remain true for the perturbational Abreu's
equation.
\end{remark}

\section {\bf  Proof of Theorem 1.1}\label{sect_3}

Let $\Omega \subset \mathbb{R}^{n}$ be a bounded convex domain with
smooth and strictly convex boundary, $\varphi $ and $\psi$
be smooth and strictly convex functions defined on an
open subset of $\mathbb{R}^{n}$ containing
$\overline {\Omega}$, satisfying
$$C_o^{-1} < \psi < C_o $$
for some constant $C_o>0$. We first consider the boundary value
problem for the perturbational Abreu equation
\begin{equation}\label{eqn_3.1}
\left\{
\begin{array}{ll}
\sum U^{ij} w_{ij} = - K ,
\;\; \det(u_{kl})= w^{-1 + \theta} & \mbox{in}\;\Omega\\
u = \varphi,\;\;\; w = \psi, & \mbox{on}\;\partial \Omega
\end{array}
\right.
\end{equation}
where $\theta$ is a very small positive number. \vskip 0.1in
\noindent We quote two results from \cite{T-W-2}.
\begin{lemma}\label{lemma_3.1}
There exists a constant $C_1>0$ such that any solution $u$ of
\eqref{eqn_3.1} satisfies
\begin{equation}\label{eqn_3.8}
C_1^{-1} \leq w \leq C_1
\end{equation}
\begin{equation}\label{eqn_3.9}
|w(\xi) - w(\xi_o)| \leq C_1 |\xi - \xi_o|\;\;\;\; \forall
\xi\;\in\; \Omega,\;\;\xi_o\;\in\;\partial\;\Omega,
\end{equation}
where $C_1$ depends only on $n$, $diam(\Omega)$, $sup_{\Omega}K$,
$sup_{\Omega}|u|$ and $C_0$.
\end{lemma}
The proof of this Lemma is the same as Lemma 7.2 in \cite{T-W-2}, we
only prove \eqref{eqn_3.8} to indicate that $C_1$ in \eqref{eqn_3.8}
is independent of $\theta$.
  Let $F=\log w $. If $F$ attains its minimum at
a boundary point, by the boundary condition \eqref{eqn_3.1} we have
$w\geq \inf_{\partial\Omega}\psi$ in $\Omega$. If $F$ attains its
minimum at an interior point $\xi \in\Omega$. By a direct
calculation, at this point we have
$$0\leq \sum u^{ij} F_{ij} \leq \frac {-K}{w^{ \theta}} \leq \frac {-k_o}{w^{ \theta}}.$$
where we used $K\geq k_o>0$. We get a contradiction. Hence $w\geq
C_o^{-1}.$

Assume $0\in \Omega.$ Let $F=\log w+\epsilon|\xi|^2,$ where
$\epsilon$ is a constant to be determined. If $F$ attains its
maximum at a boundary point, by \eqref{eqn_3.1} we have $w\leq C_o$.
If $F$ attains its maximum at an interior point $\xi_0$, Then at
this point we have
\begin{eqnarray}\label{eqn_3.4a}
&&0  = F_i=\frac {w_i}{w}+2\epsilon\xi_i,\\
\label{eqn_3.5a}&&0 \geq \sum u^{ij} F_{ij}  = \frac {-K}{w^{
\theta}} -\frac {\sum u^{ij}w_iw_j}{w^2}+2\epsilon\sum u^{ii} .
 \end{eqnarray}
Substituting \eqref{eqn_3.4a} into \eqref{eqn_3.5a} and choosing
$\epsilon=\frac{1}{5}[diam(\Omega)]^{-2}$ we have
 $$0  \geq
   \frac {-K}{w^{ \theta}} +\sum (2\epsilon -4\epsilon^2 \xi_i^2 )u^{ii}
 \geq \frac {-K}{w^{ \theta}} +\epsilon u^{ii} \geq \frac {-K}{w^{ \theta}}
 +\epsilon u^{ii}\geq  \frac {-K}{w^{ \theta}} + \epsilon w^{\frac{ 1-\theta}{n}},$$
 where we used $ \sum u^{ii}\geq  \det( u_{ij})^{-1/n}.$ Noting that $F\leq F(\xi_0),$ then
 $$w\leq e^2[C_o]^{ (n-1)\theta }\left[5\max|K|diam(\Omega)^2\right]^n.$$
 Hence as $\theta\leq 1,$ $C_1$ is independent of
 $\theta$.  \eqref{eqn_3.8} is proved. $\Box$

\begin{prop}\label{prop_3.0}
 The boundary value problem \eqref{eqn_3.1}
admits a solution $u^{(\theta)}\in C^{\infty}(\bar{\Omega})$.
\end{prop}
Please refer to Theorem 1.2 and the remark in the end of
\cite{T-W-2} for the proof. \vskip 0.1in
 Letting $\theta \rightarrow 0$, we have
\begin{prop}\label{prop_3.1}
 There is a convex function
$u\in C^{\infty}(\Omega)\bigcap C^{0}(\bar{\Omega})$ with $det
(u_{ij}) \in C^{0}(\bar{\Omega})$ such that
\begin{equation}\label{eqn_3.2}
\left\{
\begin{array}{ll}
\sum_{1\leq i,j\leq n} U^{ij} w_{ij} = - K ,
\;\; \det(u_{kl})= w^{-1} & \mbox{in}\;\Omega\\
u = \varphi,\;\;\; w = \psi, & \mbox{on}\;\partial \Omega
\end{array}
\right.
\end{equation}
\end{prop}
{\bf Proof.} By Lemma 3.1 and the comparison principle of
Monge-Ampere equations it follows that $|u^{(\theta) }|$ is
uniformly bounded. By \eqref{eqn_3.8} and a result of
Caffarelli(\cite{C1}), u is strictly convex in $\Omega$. Applying
the Caffarelli-Gutierrez theory for linearized Monge-Ampere
equations (cf. \cite{CG})
 we obtain $\det(u_{ij})\in C^{\alpha}(\Omega),$    for
some $\alpha\in(0,1)$. By Caffarelli's $C^{2,\alpha}$ estimates for
Monge-Ampere equation (\cite{C2}) we have
$$\|u\|_{C^{2,\alpha}(\Omega)}\leq C_2.$$
Following from the standard elliptic regularity theory we have
$\|u\|_{W^{4, p}(\Omega)}\leq C $. By Sobolev embedding theorem we
have
\[\|u\|_{C^{3,\alpha}(\Omega)}\leq C_2  \|u\|_{W^{4, p}(\Omega)}.
\]
Using standard bootstrap skill we conclude that
$$\|u\|_{C^{\infty}(\Omega)}\leq C_3 .$$ By
\eqref{eqn_3.9} and the bounds of $|u^{(\theta)}|$, we have $u\in
C^{0}(\bar{\Omega})$ and $\det D^2u \in C^{0}(\bar{\Omega})$. The
Proposition  is proved. $\Box$.

\begin{remark}
In \cite{CG} Caffarelli-Gutierrez  proved a H\"older estimate of
$\det(u_{ij})$ for homogeneous linearized Monge-Ampere equation
assuming that the Monge-Amp\`ere measure $\mu[u]$ satisfies some
condition, which is guaranteed by \eqref{eqn_3.8}. When $f\in
L^\infty,$ following their argument  one can obtain    H\"older
continuity of $\det(u_{ij})$.
\end{remark}

\vskip 0.1in
{\em The rest of the section is to prove  Theorem \ref{thm_1.1}. }
We consider the boundary value problem
\begin{equation}\label{eqn_3.10}
\left\{
\begin{array}{ll}
\sum U^{ij} w_{ij} = - K ,
\;\; \det(u_{kl})= w^{-1 } & \mbox{in}\;\Omega\\
u = \varphi,\;\;\; w = t, & \mbox{on}\;\partial \Omega
\end{array}
\right.
\end{equation}
 where $t>0$. We have a family of solutions
$u^{(t)}$. Let $t\rightarrow 0$. We are going to prove that:
\begin{itemize}
\item for any compact set $D\subset \Omega$ there is a subsequence
$u^{(t_i)}$, which uniformly $C^\infty$-converges on $D$; set the limit
to be $u$;
\item $u$ is smooth, strictly convex and satisfies the Abreu equation (1.2);
 \item on $\partial \Omega$
$$u = \varphi ,\;\;\;\;
|\nabla u|= \infty,\;\;\;w = 0.$$
\end{itemize}
We divide the proof into 4 steps. Each subsection consists of one step.
 \subsection{ Step One}
 We first give a
uniform estimate for $\max_{\Omega}\{u^{(t)}\} -
\min_{\Omega}\{u^{(t)}\}$. In the following calculation we will
omit the index $(t)$ to simplify notations. By adding a constant
we may assume that $\max_{\Omega}\{\varphi\} = 0.$ Then $u \leq
0$. We show that
\begin{lemma}
There is a constant $d_4>0$ independent of $t$ such that
$|u^{(t)}|\leq d_4$.
\end{lemma}
{\bf Proof. }We have
\begin{eqnarray}\label{eqn_3.3}
 - \int_{\Omega}K(u-\varphi)
&=& \int_{\Omega}\sum U^{ij} w_{ij}(u-\varphi)
= -
\int_{\Omega}\sum U^{ij}(u-\varphi)_iw_j
\\
&=& - \int_{\partial
\Omega}w\sum U^{ij}(u-\varphi)_i \gamma_j + \int_{\Omega}w\sum
U^{ij}(u-\varphi)_{ij} \nonumber\\
&\leq& - \int_{\partial \Omega}w\sum
U^{ij}(u-\varphi)_i \gamma_j + n|\Omega|,\nonumber
\end{eqnarray} where $\gamma $
denotes the unit outward normal vector on $\partial \Omega$. We
now calculate the integral $ \int_{\partial \Omega}w\sum
U^{ij}(u-\varphi)_i \gamma_j$. For any boundary point $\bar{\xi}$,
by choosing coordinates, we may assume that
$$\gamma = (0,...,0,-1).$$
 We have
$$\sum U^{ij}(u-\varphi)_i\gamma_j = (u-\varphi)_nU^{nn}.$$
Near $\bar{\xi}$ the boundary $\partial \Omega$ can be given by a
smooth and strictly convex function
$$\xi_{n} = h(\xi_1,...,\xi_{n-1})$$
with
$$h(0) = 0,\;\;\;\;\frac{\partial h}{\partial \xi_k}(0) =
0\;\;\;\;\forall k = 1,\;...,\;n-1.$$ We may write $\varphi =
u(\xi_1,...,\xi_{n-1}, h(\xi_1,...,\xi_{n-1}))$ on the boundary. By
a direct calculation we have
$$u_{kl} = \varphi_{kl} + u_n h_{kl}.$$
So
\begin{equation}\label{eqn_3.4}
U^{nn} = \det(u_{kl})|_{k,l = 1}^{n-1} \leq
a\left(\left|\frac{\partial u}{\partial
\xi_n}\right|^{n-1}+1\right),
\end{equation}
where $a$ is a constant depending on $\varphi$ and $h$.

We estimate $|\frac{\partial u}{\partial \xi_n}|$. By the same
argument  of Lemma 3.1, $\det(u_{ij})$ attains its maximum on the
boundary. Therefore there is a uniform estimate
$$\det(u_{ij}) \leq \frac{1}{t}.$$
Let $\bar{\xi}\in\partial \Omega$ be an arbitrary point. As
$\partial \Omega$ is smooth and strictly convex, by an affine
transformation ( see [L-S-C]) we may assume that $u(\bar{\xi}) =
min_{\partial \Omega}\{u\}.$ Moreover, there is point $\xi_o$ and
a constant $R
> 0$ independent of $\bar{\xi}$ such that $\Omega\subset
B_{\xi_o}(R)$, and $\Omega$ tangent to $B_{\xi_o}(R)$ at
$\bar{\xi}$. Let $\hat{u}$ be the function defined by
$$\hat{u} = \frac{1}{2} \sqrt[n]{\frac{1}{t}}(|\xi-\xi_{o}|^{2}-R^{2})+ u(\bar{\xi}).$$
We have
$$\det(u_{ij}) \leq \det (\hat{u}_{ij}),\;\;\;\; u|_{\partial \Omega} \geq \hat{u}|_{\partial \Omega}.$$
By the maximum principle, we have $\hat{u}\leq u$ in $\Omega$. As
$\hat{u}(\bar{\xi})= u(\bar{\xi})$, we get
\begin{equation}\label{eqn_3.5}
|\nabla u|\leq R\sqrt[n]{\frac{1}{t}}.
\end{equation}
It follows from \eqref{eqn_3.3}-\eqref{eqn_3.5}, we have
\begin{equation}\label{eqn_3.6}
 \int_{\partial
\Omega}w\sum U^{ij}(u-\varphi)_i \gamma_j \leq C_2
\end{equation}
for some constant $C_2>0$. Therefore
$$- \int_{\Omega}K (u-\varphi)\leq C_3$$
for some constant $C_3 >0$. Now suppose that $|u|$ attains its
maximum at $\hat\xi_o$. Let $Cone(p)$ be the cone with base
$\partial \Omega$ and vertex at $p = (\hat\xi_o, u(\hat\xi_o))$.
As $ u-\varphi \leq 0$ and $- K \leq - k_o$, we have
$$C_3 \geq - \int_{\Omega}K (u-\varphi) \geq k_o\int (\varphi - u).$$
Therefore
$$C_4 := \frac{C_3}{k_o} + \int (-\varphi)\geq \int ( - u)
\geq Vol(Cone(p)) = \frac{1}{n+1}\max_{\Omega}\{|u|\}|\Omega|.$$
Hence there is a constant $d_4$ independent of $t$ such that
$$|u^{(t)}| \leq d_4.$$
We finish the proof. $\Box$.

\v
For any compact set $D\subset \Omega$,
$$
|\nabla u^{(t)}|\leq \frac{d_4}{dist(D,\partial \Omega)}.$$ It
follows that there is a convex function $u$ defined in $\Omega$
such that for any compact set $D\subset \Omega$ there is a
subsequence $u^{(t_i)}$ converging uniformly on $D$ to $u$. We
denote the graph of $u$ by
$$M = \left\{(\xi,u(\xi))|\xi\in \Omega\right\}.$$
\subsection{ Step Two} In this
step, we prove
\begin{lemma}\label{lemma_3.2}
 For any point $\bar{p}=
(\bar{\xi}, \varphi(\bar{\xi}))$ with $\bar{\xi}\in \partial
\Omega$,
$$|\nabla u^{(t)}|(\bar{p})\rightarrow \infty$$
as $t\rightarrow 0$.
\end{lemma}
 {\bf Proof.} Suppose
that there is a point $\bar{p}= (\bar{\xi}, \varphi(\bar{\xi}))$,
where $\bar{\xi}\in \partial \Omega$, and a subsequence $t_i$ such
that $\lim _{i\rightarrow \infty}|\nabla u^{(t_i)}|(\bar{p}) = a$
for some constant $a \geq 0$. Denote $u^{(i)}= u^{(t_i)}$. $u^{(i)}$
will locally uniformly converges to a convex function $u$.
Set
$$
K_o = \frac{5 \max_{{\Omega}}\{K\}}{4d_1}.
$$
 Since
$\partial \Omega$ is smooth and strictly convex, by an affine
transformation ( see [L-S-C]) we may assume that
 \vskip 0.1in \noindent
\begin{itemize} \item
$$u^{(i)}(\bar{\xi}) = \max_{\partial \Omega}\{u^{(i)}\}
 =-\frac{n}{K_o}t_i,\;\;\;\; \bar{\xi} = (0,...,0),$$

 \item the equation of the tangent hyperplane of $\partial \Omega$ at
$\bar{\xi}$ is $\xi_1 = 0$ and  $\Omega \subset \{\xi_1 > 0\}$.
\end{itemize}

Let $\Delta(c,r)$ be a domain
defined by
$$\Delta(c,r) = \{ (\xi_1,...,\xi_n)| c(\xi_2^2 + ... + \xi_n^2) < \xi_1 < r.\}$$
Since $\partial \Omega$ is smooth, we may choose some constant $c$ and
small $r$ ( depending on $\Omega$) so  that $\Delta(c,r) \subset \Omega$.
Then $u^{(i)} < 0$ and there is a constant $b >0$ such that for
any small $\epsilon >0$
\begin{equation}\label{eqn_3.101}
|u^{(i)}| \leq b\xi_1
\end{equation}
on the part $\{\xi_1 \geq \epsilon \}\cap \Delta(c,r)$ for large $i$.  Denote the graph
of $u^{(k)}$ by
$$M^{(k)}= \left\{(\xi,u^{(k)}(\xi))|\xi\in \Omega\right\},$$

 We have
$$\sum U^{(k)ij}\left[w^{(k)} + \frac{K_o}{n}u^{(k)} \right]_{ij} = -K + K_o det(u^{(k)}_{ij})
\geq 0\quad\quad in \quad\Omega,$$
$$w^{(k)} + \frac{K_o}{n}u^{(k)} \leq 0\;\;\;on \;\;\partial
\Omega,$$ for $k$ large enough. By the maximum principle we get
$$- w^{(k)} \geq \frac{K_o}{n}u^{(k)}\quad\quad on \quad\bar{\Omega},$$
i.e.,
$$w^{(k)} \leq \frac{K_o}{n}|u^{(k)}|\quad\quad in \;\;\Omega.$$
Then
\begin{equation}\label{eqn_3.11}
 \det(u^{(k)}_{ij})\geq \frac{n}{bK_o \xi_1}
 \end{equation}
on $\{\xi_1 \geq \epsilon \}\cap\Delta(c,r)$ for $k$ large enough. We construct a
new function $\hat{u}$ as following.
We define $\hat{u}$ in $\Delta(c,r)$:
$$\hat{u} = - \left(\xi_1-c(\xi_2^2 + ... + \xi_n^2)\right)\left(-\log \xi_1\right)^{\alpha}
 + \xi_1 \left(-\log r \right)^{\alpha},$$
where $\alpha >0$ is a small positive number. It is easy to check
that
$$\hat{u}(\bar{\xi}) = 0,$$
$$\hat{u} \geq u^{(k)}\;\;\; \mbox{on}
\;\;\partial \Delta(c,r)\bigcap \{\xi_1 \geq \epsilon \},$$
\begin{equation}
\frac{\partial\hat{u}}{\partial \xi_1}(\bar{\xi}) = -(-\log \xi_1)^{\alpha} +
 (-\log r)^{\alpha}
+ \frac{\alpha}{(- \log\xi_1)^{1-\alpha}} - \frac{\alpha c(\xi_2^2 + ... +
\xi_n^2)}{(-\log \xi_1)^{1-\alpha}\xi_1}.
\end{equation}
  We now calculate
$\det(\hat{u}_{ij})$: Let $\xi = (\xi_1,\xi_2,...,\xi_n)\in
\Delta$. By taking an orthogonal transformation we may assume that
$\xi = (\xi_1,\xi_2,0,...,0)$. By a direct calculation we have, at
$\xi$,
$$\hat{u}_{11} = \frac{\alpha}{(-\log \xi_1)^{1-\alpha}\xi_1} + \alpha\frac{c\xi_2^2}{(-\log \xi_1)^{1-\alpha}\xi_1^2}
+ \frac{\alpha (1-\alpha)}{(-\log \xi_1)^{2-\alpha}\xi_1} -
\alpha(1-\alpha) \frac{c\xi_2^2}{(-\log
\xi_1)^{2-\alpha}\xi^2_1},$$
$$\hat{u}_{12} = \hat{u}_{21} = - \frac{2c\alpha \xi_2}{(-\log \xi_1)^{1- \alpha}\xi_1},\;\;\; $$
$$\hat{u}_{kk} = 2c(-\log \xi_1)^{\alpha}\;\;\;\forall \;\;k \geq 2, \;\;\;\;\hat{u}_{ij} = 0 \;\;\; \mbox{for other cases.}$$
We choose $\alpha = \frac{1}{2(n+1)}$. Obviously, $(\hat{u}_{ij})$
is positive definite for $\xi_1 < e^{-1}$, and
\begin{equation}\label{eqn_3.12}
\det(\hat{u}_{ij}) \leq \frac{\alpha (2c)^{n-1}}{(-\log \xi_1)^{1-n\alpha}\xi_1}.
\end{equation}
It is easy to see that
$$\det(u^{(k)}_{ij})\geq
\det(\hat{u}^{(k)}_{ij})\;\;\;\;\; on \;\Delta(c,r) \bigcap
\{\xi_1 \geq \epsilon \}$$ for $k$ large enough. By the maximum
principle, we have
$$\hat{u} \geq u^{(k)}\;\;\; on \;\Delta(c,r) \bigcap \{\xi_1 \geq \epsilon \}.$$
Let $k\rightarrow \infty$, $\epsilon \rightarrow 0$, we get
\begin{equation}\label{eqn_3.13}
\hat{u} \geq u \;\;\; in \;\;\Delta(c,r).
\end{equation}
From \eqref{eqn_3.101},\eqref{eqn_3.11} and \eqref{eqn_3.13} we get a contradiction. The lemma
is proved. $\Box$
\vskip 0.1in \noindent Using the same method
we can also prove
\begin{lemma}\label{lemma_3.3}
 For any point $\xi\in \Omega$ and any support hyperplane $H$ of $M$
 at $p=(\xi,u(\xi))$,
 $\mathrm{dist}( H,\partial M) >0.$
 \end{lemma}

Denote by $f^{(t)}(x)$ the Legendre transformation functions of
$u^{(t)}$. For any large $R>0$, by Lemma \ref{lemma_3.2}, there is a
$t_{0}>0$ such that $f^{(t)}$ is defined on the disk $B_0(R)$ for
$0< t < t_{0}$. Since $|\nabla f^{(t)}|$ are uniformly bounded by
$diam (\Omega)$, there is a subsequence $f^{(t_i)}$ converging to
a convex function $f$ defined on $B_0(R)$. Let $R\rightarrow
\infty$. By choosing subsequences we conclude that there exists a
convex function $f$ defines on whole $\mathbb{R}^n$ such that $f^{(t)}$
locally uniformly converges to $f$.

\subsection{Step Three}
We first prove
\begin{lemma}\label{lemma_3.4} We assert that
\begin{itemize} \item For any point $\xi \in \Omega$ and any
support hyperplane $H$ of $M$ at $p=(\xi,u(\xi))$
$$\dim (M\bigcap H)\leq n-1;$$
\item For any point $\xi \in \Omega$ and any ball
$B_{\xi}(\delta')\subset \Omega$ with radius $\delta'$ around $\xi$, there exists a point
$\xi_o\in B_{\xi}(\delta')$ such that $u$ has second derivatives
and strictly convex at $\xi_o$.
\end{itemize}
\end{lemma}
{\bf Proof. } The first claim follows from
the weakly continuous of Monge -Ampere measure. We prove the
second one. Since $u$ is convex, it has second order derivatives
almost everywhere. Let $G \subset B_{\xi}(\delta')$ be the set
where $u$ has second order derivatives. Then $| B_{\xi}(\delta')
-G| = 0$ . Let $O$ be an open subset of $B_\xi(\delta')$
 such that $ B_{\xi}(\delta')-G
\subset O$ with $|O|\leq \epsilon'$. We choose $\epsilon'$ so
small that $|B_{\xi}(\delta')- O| >
\frac{1}{2}|B_{\xi}(\delta')|$. By  the weak convergence
of the Monge-Ampere measure we have
\begin{equation}\label{eqn_3.14}
\int_{B_{\xi}(\delta')- O}\det(u_{kl}) > \frac{1}{2}d_1
|B_{\xi}(\delta')|.
\end{equation}
 We claim that there exists a
point $\xi^o \in B_{\xi}(\delta')- O$ such that
$$\det(u_{kl})(\xi^o) \geq
\frac{d_1|B_{\xi}(\delta')|}{2|B_{\xi}(\delta')- O|}.$$
Otherwise,
we would have
$$
\int_{B_{\xi}(\delta')- O}\det(u_{kl}) \leq \frac{1}{2} d_1
|B_{\xi}(\delta')|,$$
which contradicts to \eqref{eqn_3.14}. As
$u$ has second order derivatives at $\xi_o$, there is a constant
$B$ such that $\lambda_k \leq B$ for all $k$, where $\lambda_k$
denote the eigenvalues of $u_{kl}(\xi_o)$. Denote by $\lambda_1$
the least eigenvalue of $u_{kl}$. Then
$$\frac{d_1}{4} < \frac{d_1|B_{\xi}(\delta')|}{2|B_{\xi}(\delta')- O|}<  \lambda_1 B^{n-1},$$
which implies that $u$ is strictly convex at $p_o$. $\Box$.

\vskip
0.1in \noindent
Without loss of generality we assume that $u$ has
second order derivatives and strictly convex at $0$, and
$$u \geq 0,\;\;\; u(0) = 0.$$
 Then there is a positive
number $a
> 0$ such that
\begin{itemize}
\item $\partial M\bigcap \{\xi_{n+1} = a\} \ne \phi$ \item $\partial M\bigcap
\{\xi_{n+1} = a-\epsilon\} = \phi$ for any $\epsilon
>0$.
\end{itemize}
\begin{lemma}\label{lemma_3.5}
 $u$ is smooth and strictly
convex in $S_u(0,a)$.
\end{lemma}
{\bf Proof.} Suppose that $u^{(i)}$ locally uniformly converges
to $u$. Since $u$ is strictly convex at $0$, there is a small
positive number $\epsilon''$ and $b' >0$ such that
$\bar{S}_{u^{(i)}}(0,\epsilon'')$ is compact and
$$
\sum \left(\frac{\partial u^{(i)}}{\partial \xi_k}\right)^2\leq b'
\;\;\; in\;\;\; S_{u^{(i)}}(0,\epsilon'')$$ for large $i$. By
Lemma \ref{lemma_2.1} and Lemma \ref{lemma_2.2} we have uniform
estimates for $\det(u^{(i)}_{kl})$ on $S_{
u^{(i)}}(0,\frac{1}{2}\epsilon'')$ both from above and below. Then
we use the Caffarelli-Gutierrez theory and the Caffarelli-Schauder
estimate to conclude that $\{u^{(i)}\}$
 smoothly converges $u$. Therefore,
$u$ is  a smooth and
strictly convex function in $S_u(0,\frac{1}{2}\epsilon'')$ and $u$
satisfies the Abreu equation $\s(u) = K$.

Let $f^{(i)}(x)$ be the legendre transformations of $u^{(i)}$.
 Then $f^{(i)}$ locally
uniformly converge to a convex function $f(x)$ defined on the
whole $\mathbb{R}^n$. Furthermore, in a neighborhood of $0$,
$f(x)$ is a smooth strictly convex function such that its
Legrendre transform $u$ satisfies the Abreu equation. Denote
$$\widetilde M ^{(\infty)} = \left\{(x_1,...,x_n,f(x))\right\}.$$
By the convexity of $f^{(i)}$ there is a constant $b''>0$ such
that
$$\frac{\sum x_k^2}{(1+ f^{(i)})^2} \leq b$$
for any $i$.

Let $\Gamma\subset S_u(0,a)$ be the set of points where
$u$ is smooth and strictly convex. By the same argument given above,
we know that
$u^{(i)}$ smoothly converges to $u$ on $\Gamma$ and $\Gamma$
is open. If $\Gamma=S_u(0,a)$, the lemma is proved.
Otherwise, there is a constant $c$ with $0<c<a$ such that
$S_u(0,c)\subset \Gamma$ but $\bar S_u(0,c)$ is not.
When this is the case,
we choose a point $\bar \xi\in \bar S_u(0,c)\setminus \Gamma$.
we may choose a point $\xi_o \in S_u(0,c)$ close
to $\bar{\xi}$ and a small $\delta>0$ such that
$$
\bar{\xi} \in S_{u}(\xi_o,\delta)\subset S_u(0,a).
$$
It is then again that we repeat the above argument to show
that
$u$ is smooth and strictly
convex in the section $S_{u}(\xi_o,\delta)$:
let $\hat u^{(i)}$ be the normalization of $u^{(i)}$ at $\xi_o$;
let $\hat x=\partial \hat u^{(i)}/\partial \xi$ and $\hat f^{(i)}$
be the Legendre transform of $\hat u^{(i)}$;
then by Lemma \ref{lemma_2.35}
$$
\frac{|\hat x|^2}{(d+ \hat f^{(i)})^2} \leq b.
$$
Then applying Lemma \ref{lemma_2.3} we get uniform bounds of $\det
(u^{(i)}_{ij})$ at $S_{\hat
u^{(i)}}(\xi_o,\delta'),\delta'<\delta$ for large $i$; therefore
$\hat u^{(i)}$, and so $u^{(i)}$, smoothly converges in this
domain. This contradicts to the assumption that $\bar\xi\notin
\Gamma$. $\Box$

\vskip 0.1in \noindent We are going to prove that $u$ is smooth
and strict convex in $\Omega$. If not, then there is a open set
$U\subset \bar{U}\subset \Omega$ such that
\begin{itemize}
\item $u$ is smooth and strictly convex in $U$; \item there is a
point $\bar{q}=(\bar{\xi}, u(\bar{\xi}))$, $\bar{\xi}\in \bar{U}$,
and a line segment $L$ such that $\bar{p}\in L \subset M$.
\end{itemize}
By Lemma 3.5 for any support hyperplane $P$ of $M$ at $\bar{q}$ we
have $\mathrm{dist}(\partial M, H) >0$. We can choose a point $\xi_o \in U$, very
close to $\bar{\xi}$, and a positive number $a_1$ such that
$$\bar{\xi}\in S_{u}(\xi_o,P, \frac{1}{2}a_1),
\mbox{ and, } \partial M \bigcap \{\xi_{n+1}= u(\xi_o) + P(\xi - \xi_o) + b\} = \phi$$
for any $b\leq a_1$, where $(P,-1)$ is the normal vector of $M$ at
$(\xi_o,u(\xi_o))$ and
$$S_u(\xi_o,P,\frac{1}{2}a_1):=\{\xi\in \Omega|u(\xi) < u(\xi_o) + P(\xi - \xi_o) +
\frac{1}{2}a_1\}.$$
 By Lemma \ref{lemma_3.5}, we know that  $u$ is
smooth and strictly convex in $S_u(\xi_o, a_1)$, we get a
contradiction. It follows that $u$ is smooth and strictly convex
in $\Omega$.

\subsection{Step Four}
Since $u$ is smooth and strictly convex in $\Omega$, it can be extended over $\partial\Omega$.
Let $\varphi'$ be the function. Then obviously, on $\partial\Omega$, $\varphi'\leq \varphi$.
We prove that
\begin{lemma}
$\varphi'=\varphi$.
\end{lemma}
{\bf Proof. }
Suppose that there is a point $\bar{\xi}\in \partial
\Omega $ such that $\varphi'(\bar{\xi})< \varphi(\bar{\xi})$.
Without loss of generality, we assume that
\begin{itemize} \item
$\varphi(\bar{\xi}) = 1,\;\;\; \varphi'(\bar{\xi})= 0,$
 \item  $\bar{\xi} = (0,...,0),$ and the equation of the tangent
 hyperplane of $\partial \Omega$ at
$\bar{\xi}$ is $\xi_1 = 0$ and  $\Omega \subset \{\xi_1 > 0\}$,
\item restricting to $\{\xi_1 \leq \epsilon\}\bigcap \partial
\Omega$, $\varphi' < \frac{1}{10}$.
\end{itemize}
We construct a function $\tilde{u}$
$$\tilde{u} = 2u - b\xi_1 + \frac{1}{2}.$$
By choosing $b$ large, we have $\tilde{u} + \frac{1}{3} \leq u $ on
$\partial \Delta'$, where
$$\Omega' := \{\xi \in \Omega | \xi_1\leq \epsilon\}.$$
For any positive number $\delta >0$, let $D_{\delta} = \{\xi\in
\Omega'| \mathrm{dist}(\xi, \partial \Omega') \geq \delta \}$. Then for
$\delta$ small enough and $ i$ large enough we have
$$\tilde{u} < u\;\;\;on \;\;\partial D_{\delta},$$
$$\det(\tilde{u}_{kl}) >\det(u^{(i)}_{kl}).$$
It follows that $u^{(i)} \geq \tilde{u}$ on $\overline
{D}_{\delta}$. Letting $\delta \rightarrow 0$, $i\rightarrow
\infty$ we get $\tilde{u}\leq u$ in $\Delta'$. But
$\tilde{u}(\bar{\xi}) = \frac{1}{2} > \varphi'(\bar{\xi})$, and
both $\tilde{u}$ and $u$ is smooth in the interior of $\Delta'$,
we get a contradiction. So $\varphi' = \varphi $ on $\partial
\Omega$. The claim $w = 0$ on $\partial \Omega$ follows from Lemma
\ref{lemma_2.7}. The claim $|\nabla u| = \infty$ follows from the fact that
$f$ is defined on the whole $\mathbb{R}^n$.

\end{document}